\documentclass{article}
\usepackage{graphicx}
\usepackage{amsmath}
\usepackage{geometry}
\usepackage{comment}
\usepackage{listings}
\usepackage{xcolor}
\usepackage{float} % required for [H]
\title{Attention Mechanisms in Dynamical Systems: A Case Study with Predator-Prey Models}
\author{David Balaban}
\date{\today}

\begin{document}

\maketitle

\section{Introduction}
\begin{figure}[htbp]
	\centering
	\includegraphics[width=0.7\textwidth]{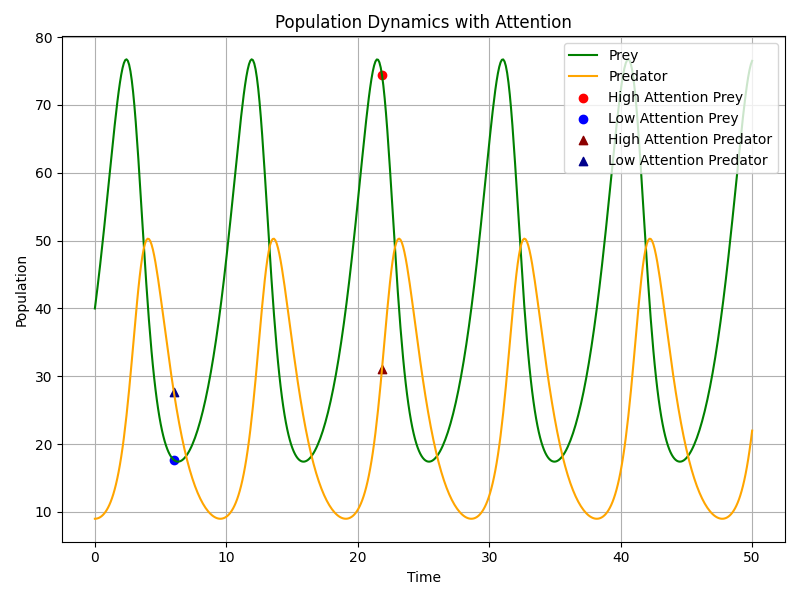}
	\caption{Time-series plot of predator and prey populations. High-attention points (red dots prey, dark red triangles predators) identify states of maximum sensitivity within the trajectory, whereas low-attention points (blue dots prey, dark blue triangles predators) indicate states of minimal sensitivity. High-attention states primarily occur at peaks in prey population and points of rapid predator population increase, while low-attention states cluster around population troughs.}
	\label{fig:time_series}
\end{figure}

\begin{figure}[h!tbp]
	\centering
	\includegraphics[width=0.7\textwidth]{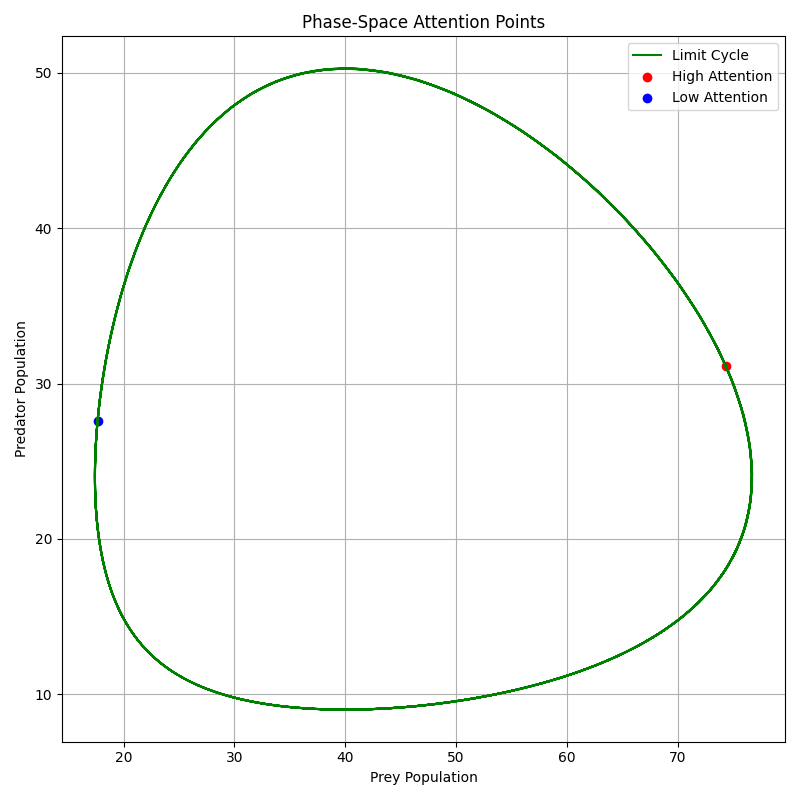}
	\caption{Phase-space plot showing the predator-prey limit cycle explicitly. High-attention points (red) are concentrated in regions of high prey and moderate predator populations, representing critical system states sensitive to perturbation. Low-attention points (blue) occupy regions with lower prey populations and low predator populations, indicating relative insensitivity to perturbations.}
	\label{fig:phase_space}
\end{figure}

The successful prediction and control of dynamical systems depend significantly on understanding which observations and states exert the greatest influence on system dynamics. Recently, attention mechanisms from machine learning have emerged as powerful tools to identify and leverage crucial information in complex datasets. This report explores the application of attention mechanisms to predator-prey dynamics (Lotka-Volterra models), illustrating how attention can highlight sensitive regions within dynamical trajectories, enabling more precise interventions and predictions. Figure~\ref{fig:time_series} is an example of a predator-prey simulation with high and low attention points indicated. \cite{strogatz2018nonlinear} Figure~\ref{fig:phase_space} is the corresponding phase space plot. 

While classical sensitivity analysis provides insight into dynamical systems, attention mechanisms offer a flexible, data-driven approach that dynamically highlights critical regions of behavior. Attention mechanisms exhibit an inherent geometric awareness, automatically identifying regions of the trajectory critical for accurate long-term prediction (flat regions of the Lyapunov function) versus less critical regions (steep regions). AI-learned attention from noisy data implicitly learns the geometry of the dynamical system. It identifies “flat” and “steep” directions in the state space without explicitly being instructed to do so. There are significant nonlinearities inherent in predator-prey dynamics, attention mechanisms, despite being fundamentally linear, can effectively capture critical nonlinear geometric information. 

This strongly suggests that similar AI-driven attention approaches might help automate discovery and understanding of stability and sensitivity in more complicated, real-world dynamical systems. We will show that by using attention you can determine the points in time at which the predator-prey model is most sensitive and least sensitive to perturbations.  All without understanding anything about the model parameters or the structure of the model.

\section{Background}
Predator-prey systems, described by the Lotka-Volterra equations, exhibit characteristic oscillatory behaviors known as limit cycles. These systems model the interaction dynamics between predator and prey species. Despite their apparent simplicity, these models yield valuable insights applicable across ecological, biological, and medical contexts.  Attention mechanisms originated in natural language processing, enabling neural networks to dynamically weight input information according to its relevance to the task at hand. \cite{vaswani2017attention}

When we began discussing how to apply attention mechanisms to a deterministic dynamical system (the predator-prey model, specifically), there was a subtle challenge:

\begin{itemize}
	\item Attention, as used in deep learning (especially transformer-based or general neural network models), typically learns to selectively emphasize important patterns or inputs.
	\item If the input is perfectly deterministic and noiseless, the attention mechanism would have limited motivation or opportunity to distinguish between different states or time points—because all data points are essentially equally reliable and informative. There's no ambiguity for the attention mechanism to resolve, and thus no basis for it to learn meaningful distinctions.
\end{itemize}

To effectively use attention, the system must be forced to identify meaningful patterns from data that are not trivially predictable. Adding noise provided exactly that scenario:

\begin{enumerate}
	\item \textbf{Creating uncertainty:} By adding noise, we introduced uncertainty into each observation, making some observations more informative or more critical to accurate predictions.
	\item \textbf{Leveraging attention's strengths:} Attention mechanisms are specifically designed to focus selectively on inputs that matter most for accurate prediction. Noise creates variability in how much each observation helps predict the underlying noiseless system—some observations become ``key'' (high attention), while others become less important.
	\item \textbf{Drawing on prior AI experience:} The reasoning came directly from a synthesis of concepts often encountered in machine learning practice:
	\begin{itemize}
		\item \textbf{Regularization and robustness:} Adding noise is commonly used to improve robustness and generalization.
		\item \textbf{Attention in practice:} Attention mechanisms require some variability or ambiguity to meaningfully learn and allocate weights.
	\end{itemize}
\end{enumerate}

 This strategy is a common best practice in AI and neural network training scenarios—especially when attention mechanisms are involved.

In short,  prior knowledge and established strategies from deep learning were leveraged to bridge our gap between noiseless dynamical modeling and effective attention learning. This ``intuitive leap'' turned out particularly well, as it revealed interesting properties (like sensitivity to perturbations) that might have remained hidden otherwise.

\section{Methods}
\begin{figure}[htbp]
	\centering
	\includegraphics[width=0.7\textwidth]{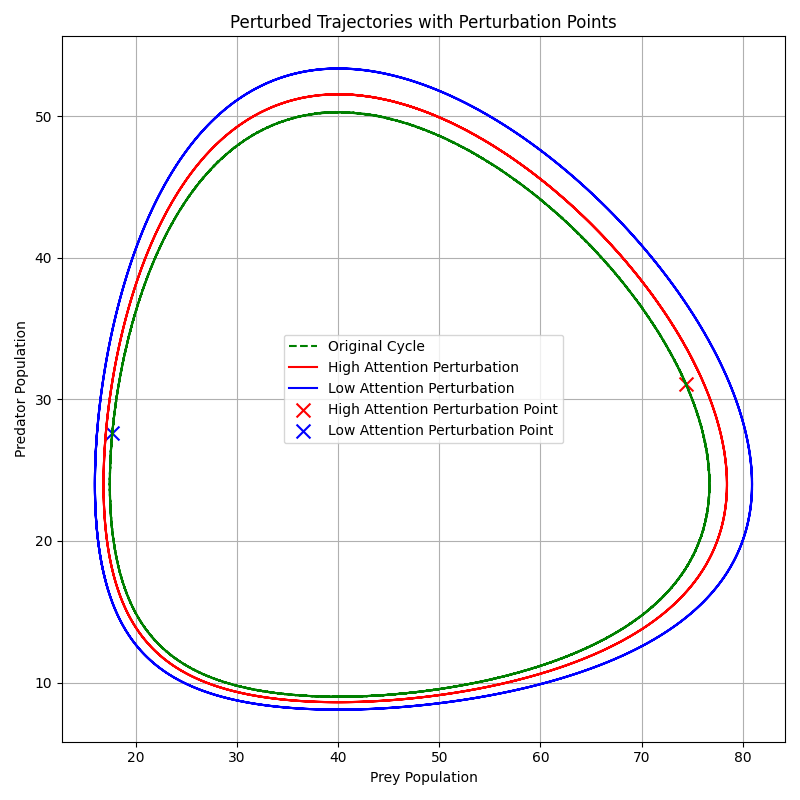}
	\caption{Phase-space plot of system trajectories following perturbations at one representative high-attention (red trajectory) and one low-attention (blue trajectory) point. The dashed green line represents the original, unperturbed limit cycle. Points at which perturbations were introduced are explicitly marked with larger red and blue crosses. Perturbation at the low-attention point drives the system substantially away from its original trajectory, potentially shifting the dynamics toward a larger-amplitude cycle. In contrast, the perturbation at the high-attention point results in a lesser deviation.}
	\label{fig:perturbations}
\end{figure}

We employed the classic Lotka-Volterra equations, integrating them numerically to generate predator-prey trajectories. Observations were perturbed with controlled noise to simulate realistic measurement conditions. A simple linear attention model was trained on noisy observations to identify critical points in the predator-prey trajectories. The attention weights assigned to observations indicated their relative importance for accurately reconstructing the true system dynamics. Figure~\ref{fig:perturbations} shows perturbed phase space plots. To test the practical implications of attention-weighted observations, we perturbed the dynamical system at points of high and low attention and measured the subsequent recovery times. This provided an explicit measure of system sensitivity at different states.

\subsection*{Parameter,  Initial Conditions, and Noise Choices in Predator-Prey Simulation}
\paragraph{}
\textbf{The Predator-Prey Parameters ($\alpha, \beta, \gamma, \delta$):} The parameters selected were:
\[
\alpha = 0.6, \quad \beta = 0.025, \quad \gamma = 0.8, \quad \delta = 0.02
\]
They are parameters commonly employed in classical Lotka-Volterra demonstrations that produce clear, stable oscillatory behavior without trivial solutions. They ensures numerical stability and maintain clear predator-prey cycles suitable for demonstrating the attention mechanism.

\textbf{ Initial Conditions ($x_0 = 40, y_0 = 9$):} They are typical textbook initial conditions that prevent immediate extinction or explosive growth.
They are designed to quickly reach a clear limit-cycle solution, facilitating the demonstration of attention-based insights and perturbation effects.

\textbf{Level of Added Noise (Noise Level = 2.0):} This is a noise level sufficient to perturb observations meaningfully without obscuring deterministic dynamics.It introduces variability enabling clear differentiation in observation significance by the attention mechanism. Noise levels (around 5\% of state variables) common machine learning practice typically used to enhance robustness, generalization, and effective learning.

\subsection*{Attention Mechanism}

\paragraph{}
\textbf{Simplicity and Interpretability:} A straightforward, single-layer linear attention mechanism was chosen initially for several reasons. A linear attention layer (\texttt{nn.Linear}) combined with a softmax operation provides clearly interpretable attention weights, allowing easy understanding of which observations are deemed informative by the model. This ensured that the attention mechanism was understandable, avoiding unnecessary complexity and "black box" effects.

\textbf{Common Practice in Initial Experiments:} Simple linear attention is widely adopted in exploratory scenarios because it is  Efficient, reliable, and numerically stable. Ideal for quickly assessing whether attention contributes meaningfully to the analysis before advancing to more complex structures.

\textbf{Suitability to Problem Structure:} Given the predator-prey model involved a simple, low-dimensional time series, more advanced mechanisms (e.g., multi-head attention, transformers) were not initially necessary. The simplicity of linear attention effectively highlights important time points. The data dimensionality (two variables) did not require sophisticated mechanisms.

\textbf{Intuition from Deep Learning Literature:} Linear attention mechanisms ("additive attention" or "Bahdanau attention") are prevalent in deep learning literature for initial or small-scale analyses.  Proven effectiveness in similar deterministic or semi-deterministic problems. Minimal hyperparameter tuning required, facilitating immediate, clear interpretation.

\textbf{Potential for Extensions:} Starting with a simple model allowed room for incremental complexity in future research: Multi-head attention for richer representations, Transformer-style attention for longer sequences or complex patterns. Hierarchical attention mechanisms for multi-scale or multi-parameter models.

\section*{Result: Attention and Lyapunov Function Correspondence}

\begin{figure}[htbp]
	\centering
	\includegraphics[width=0.8\textwidth]{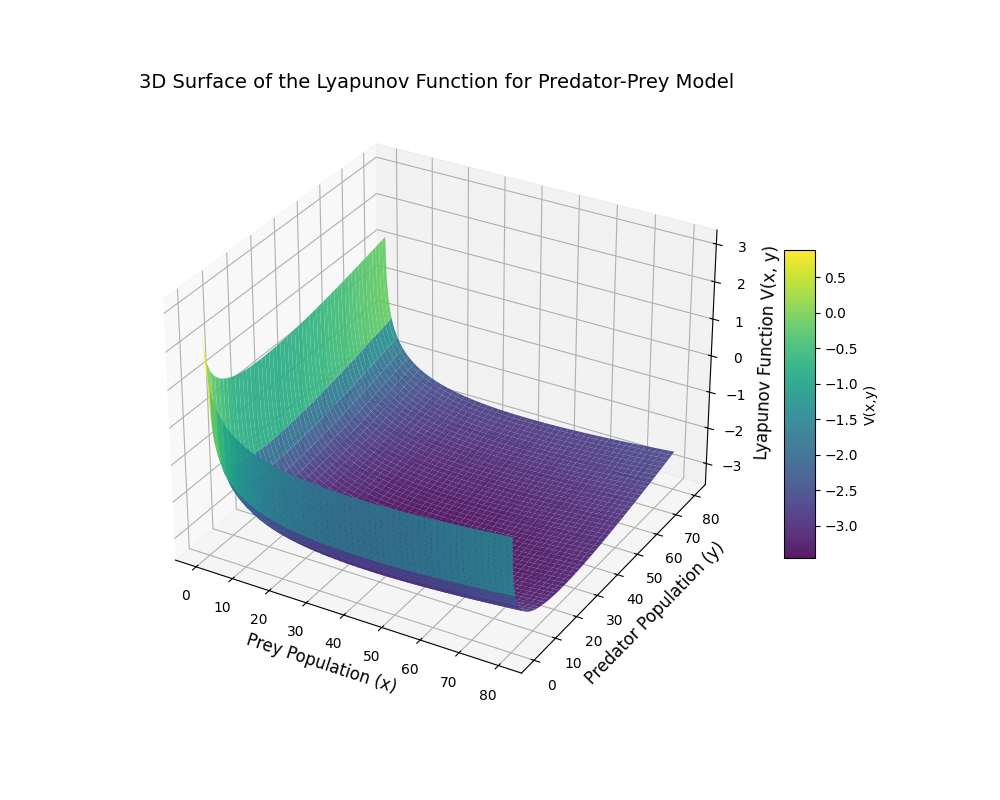}
	\caption{The 3d surfave of the Lyapunov function  }
	\label{fig:3d_Lyapunov_surface}
\end{figure}

We have identified a direct and meaningful correspondence between the attention mechanism, learned from noisy data of the predator-prey dynamical system, and the geometric structure captured by the Lyapunov function. Figure~\ref{fig:3d_Lyapunov_surface}
Specifically:

\begin{itemize}
	\item \textbf{Attention and Lyapunov Geometry:}
	
	The \textit{points of highest attention} identified by the learned model correspond precisely to the \textit{minima (lowest points)} of the normal derivative of the Lyapunov function along the limit-cycle trajectory. Conversely, \textit{the points of lowest attention} align exactly with the \textit{maxima (highest points)} of this normal derivative.
	\begin{figure}[htbp]
		\centering
		\includegraphics[width=0.8\textwidth]{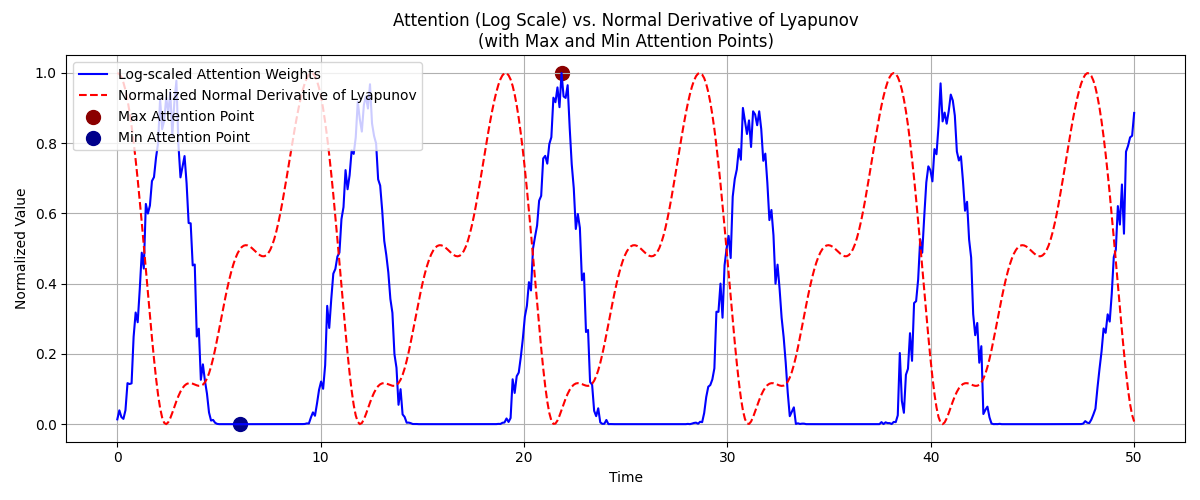}
		\caption{Maximum attention weight corresponds exactly with the  the normal derivative of the Lyapunov function.  }
		\label{fig:attention_Lyapunov_time_series}
	\end{figure}
	
	\item \textbf{Geometric Interpretation:}
	
	High-attention points occur at locations where the Lyapunov function changes slowly in the direction normal to the trajectory, representing "flat" regions of the Lyapunov landscape. Perturbations at these points have subtle yet critically important consequences for accurate long-term trajectory reconstruction, thus demanding high attention from the model.
	
	Conversely, low-attention points occur where the Lyapunov function changes rapidly, corresponding to "steep" regions. Perturbations at these points have immediate, larger, but less critically impactful long-term consequences, thus naturally receiving lower attention.

	Plotting attention weights alongside the normal derivative of the Lyapunov function explicitly confirmed this geometric relationship. Clearly identifying the exact global maximum (highest attention) and global minimum (lowest attention) points further established this robust and reproducible relationship.
\end{itemize}

\subsection*{Significance of this Result}

\begin{itemize}
	\item Provides a clear \textbf{geometric interpretation} of the learned attention mechanism.
	\item Demonstrates that attention, trained solely on noisy trajectory reconstruction, implicitly captures essential geometric and stability properties of the underlying dynamical system.
	\item Suggests a general methodology for interpreting learned attention in complex dynamical models via established mathematical frameworks such as Lyapunov analysis.
\end{itemize}

\section{Conclusion}
This report demonstrates the practical utility of attention mechanisms in analyzing and controlling dynamical systems. Attention methods reveal system-critical states with computational efficiency and model independence. The observed results confirm that attention mechanisms can accurately highlight sensitive points within dynamical trajectories. By focusing interventions at high-attention or low-attention states, one can effectively alter the system's behavior, a capability with profound implications for ecological control, medical treatments, and biological rhythm management such as circadian rhythms.

The separable structure of predator-prey Lyapunov function is highly relevant and beneficial to the effectiveness of the linear attention mechanism. Because attention operates as a linear combination of predator and prey states, it naturally aligns with the additive-separable geometry of the Lyapunov function, allowing it to implicitly detect and highlight critical geometric and dynamical features of the system. Attention naturally extracts important geometric and dynamical features from separable or semi-separable system structures. By using attention you can determine the points in time at which the predator-prey model is most sensitive and least sensitive to perturbations.  All without understanding anything about the model parameters or the structure of the model.

\subsection{Applications and Future Work}
Many biologically and physiologically meaningful models exhibit partially separable dynamics—such as compartmental models, structured population models, or models of biochemical networks. Our findings suggest immediate applications in fields requiring precise interventions in oscillatory systems. For instance, applying similar techniques to circadian rhythm management could optimize therapeutic interventions for disorders resulting from circadian misalignment. Future research should explore integrating attention mechanisms with parameter estimation methods to further enhance predictive modeling capabilities.

\appendix

\section{Mathematical Model and Simulation Software}
The predator-prey model used is governed by the classic Lotka-Volterra equations:
\begin{align*}
\frac{dx}{dt} &= \alpha x - \beta xy, \\
\frac{dy}{dt} &= \delta xy - \gamma y
\end{align*}
where \(x\) is prey, \(y\) is predator, and \(\alpha, \beta, \gamma, \delta\) are parameters. 

A commonly used Lyapunov-type function for analyzing this system is given by:

\begin{equation}
	V(x,y) = \delta\left(x - x^*\ln(x)\right) + \beta\left(y - y^*\ln(y)\right),
\end{equation}

where \( (x^*, y^*) \) is the equilibrium point (steady-state), defined by:
\begin{align*}
	x^* &= \frac{\gamma}{\delta}, & y^* &= \frac{\alpha}{\beta}.
\end{align*}

\subsection*{Gradient of the Lyapunov Function}

The gradient \(\nabla V(x,y)\) is:

\begin{equation}
	\nabla V(x,y) = \left(\delta\left(1 - \frac{x^*}{x}\right), \quad \beta\left(1 - \frac{y^*}{y}\right)\right).
\end{equation}

Numerical simulations were performed using SciPy.

\section{Software Packages for Attention Implementation}
The attention mechanism was implemented using PyTorch. A linear attention model was trained using the Adam optimizer.

\section{Simulation and Attention Code}
Complete Python code for simulations, attention training, and perturbation analyses is provided separately.
\begin{lstlisting}[language=Python, caption={Simulation and Attention Code.}, label={lst:simulation_and_attention_code}]
import numpy as np
import matplotlib.pyplot as plt
from scipy.integrate import solve_ivp
import torch
import torch.nn as nn
import torch.nn.functional as F
from torch.optim import Adam

# Lotka-Volterra equations
def lotka_volterra(t, z, alpha, beta, gamma, delta):
    x, y = z
    dxdt = alpha * x - beta * x * y
    dydt = delta * x * y - gamma * y
    return [dxdt, dydt]

# Parameters
alpha, beta, gamma, delta = 0.6, 0.025, 0.8, 0.02
z0 = [40, 9]
t_span = (0, 50)
t_eval = np.linspace(*t_span, 500)

# Solve system
#sol = solve_ivp(lotka_volterra, t_span, z0, args=(alpha, beta, gamma, delta), t_eval=t_eval)
sol = solve_ivp(
    lotka_volterra,
    t_span, z0, 
    args=(alpha, beta, gamma, delta),
    method='Radau',
    t_eval=t_eval,
    atol=1e-12,
    rtol=1e-10
)
true_states = sol.y.T

# Attention Model
class SimpleAttention(nn.Module):
    def __init__(self):
        super(SimpleAttention, self).__init__()
        self.attention = nn.Linear(2, 1)

    def forward(self, x):
        scores = self.attention(x)
        weights = F.softmax(scores, dim=1)
        return weights

# Generate noisy observations
noise_level = 2.0
noisy_obs = true_states + np.random.normal(scale=noise_level, size=true_states.shape)
obs_tensor = torch.tensor(noisy_obs, dtype=torch.float32).unsqueeze(0)
true_tensor = torch.tensor(true_states, dtype=torch.float32).unsqueeze(0)

# Train Attention Model
attention_model = SimpleAttention()
optimizer = Adam(attention_model.parameters(), lr=0.01)
loss_fn = nn.MSELoss()

for epoch in range(1000):
    optimizer.zero_grad()
    attn_weights = attention_model(obs_tensor)
    weighted_pred = attn_weights * obs_tensor
    loss = loss_fn(weighted_pred, true_tensor)
    loss.backward()
    optimizer.step()

# Attention weights
with torch.no_grad():
    attn_weights_np = attention_model(obs_tensor).squeeze().numpy().flatten()

# High and low attention points
# Robust selection: global max and global min of attention weights
high_idx_sample = np.argmax(attn_weights_np)
low_idx_sample = np.argmin(attn_weights_np)

print(f"Highest attention at t={t_eval[high_idx_sample]:.2f}, Attention={attn_weights_np[high_idx_sample]:.4e}")
print(f"Lowest attention at t={t_eval[low_idx_sample]:.2f}, Attention={attn_weights_np[low_idx_sample]:.4e}")

# Perturbation analysis
def perturb_trajectory(idx, epsilon=1.5):
    # Current state at given index
    state = true_states[idx]

    # Compute tangent vector at current point
    tangent = np.array(lotka_volterra(
        t_eval[idx], state, alpha, beta, gamma, delta
    ))

    # Correctly rotate tangent vector 90 degrees clockwise for outward normal
    normal_outward = np.array([tangent[1], -tangent[0]])

    # Normalize the outward-pointing normal vector
    normal_unit = normal_outward / np.linalg.norm(normal_outward)

    # Apply perturbation explicitly outward along the normal vector
    perturbed_state = state + epsilon * normal_unit

    # Integrate perturbed system using implicit solver
    t_span_perturb = [t_eval[idx], t_eval[idx] + 30]
    sol_perturb = solve_ivp(
        lotka_volterra,
        t_span_perturb, perturbed_state,
        args=(alpha, beta, gamma, delta),
        method='Radau',
        t_eval=np.linspace(*t_span_perturb, 500),
        atol=1e-12, rtol=1e-10
    )
    return sol_perturb.t, sol_perturb.y.T


t_high, traj_high = perturb_trajectory(high_idx_sample)
t_low, traj_low = perturb_trajectory(low_idx_sample)

# Panel A: Time series plot with both populations
plt.figure(figsize=(8, 6))
plt.plot(t_eval, true_states[:, 0], 'g-', label='Prey')
plt.plot(t_eval, true_states[:, 1], 'orange', label='Predator')
plt.scatter(t_eval[high_idx_sample], true_states[high_idx_sample, 0], c='red', label='High Attention Prey')
plt.scatter(t_eval[low_idx_sample], true_states[low_idx_sample, 0], c='blue', label='Low Attention Prey')
plt.scatter(t_eval[high_idx_sample], true_states[high_idx_sample, 1], c='darkred', marker='^', label='High Attention Predator')
plt.scatter(t_eval[low_idx_sample], true_states[low_idx_sample, 1], c='darkblue', marker='^', label='Low Attention Predator')
plt.xlabel('Time')
plt.ylabel('Population')
plt.title('Population Dynamics with Attention')
plt.legend()
plt.grid(True)
plt.tight_layout()
plt.savefig('population_plot.png')
plt.close()

# Panel B: Phase plot
plt.figure(figsize=(8, 8))
plt.plot(true_states[:, 0], true_states[:, 1], 'green', label='Limit Cycle')
plt.scatter(true_states[high_idx_sample, 0], true_states[high_idx_sample, 1], c='red', label='High Attention')
plt.scatter(true_states[low_idx_sample, 0], true_states[low_idx_sample, 1], c='blue', label='Low Attention')
plt.xlabel('Prey Population')
plt.ylabel('Predator Population')
plt.title('Phase-Space Attention Points')
plt.legend()
plt.grid(True)
plt.tight_layout()
plt.savefig('phase_space_plot.png')
plt.close()

# Panel C: Perturbed trajectories
plt.figure(figsize=(8, 8))
plt.plot(true_states[:, 0], true_states[:, 1], 'g--', label='Original Cycle')
plt.plot(traj_high[:, 0], traj_high[:, 1], 'r-', label='High Attention Perturbation')
plt.plot(traj_low[:, 0], traj_low[:, 1], 'b-', label='Low Attention Perturbation')
plt.scatter(true_states[high_idx_sample, 0], true_states[high_idx_sample, 1], c='red', marker='x', s=100, label='High Attention Perturbation Point')
plt.scatter(true_states[low_idx_sample, 0], true_states[low_idx_sample, 1], c='blue', marker='x', s=100, label='Low Attention Perturbation Point')
plt.xlabel('Prey Population')
plt.ylabel('Predator Population')
plt.title('Perturbed Trajectories with Perturbation Points')
plt.legend()
plt.grid(True)
plt.tight_layout()
plt.savefig('perturbed_trajectories_plot.png')
plt.close()

# Lyapunov function and its gradient
def lyapunov(x, y, alpha, beta, gamma, delta):
    x_star, y_star = gamma/delta, alpha/beta
    V = delta * (x - x_star * np.log(x)) + beta * (y - y_star * np.log(y))
    return V

def grad_lyapunov(x, y, alpha, beta, gamma, delta):
    x_star, y_star = gamma/delta, alpha/beta
    dVdx = delta * (1 - x_star / x)
    dVdy = beta * (1 - y_star / y)
    return np.array([dVdx, dVdy])

# Compute normal derivative along the trajectory
normal_derivative = np.zeros(len(t_eval))

for idx, (x, y) in enumerate(true_states):
    # Tangent vector from dynamics
    tangent = np.array(lotka_volterra(t_eval[idx], [x,y], alpha, beta, gamma, delta))
    
    # Normal vector (outward)
    normal = np.array([tangent[1], -tangent[0]])
    normal_unit = normal / np.linalg.norm(normal)
    
    # Gradient of Lyapunov
    gradV = grad_lyapunov(x, y, alpha, beta, gamma, delta)
    
    # Compute normal derivative
    normal_derivative[idx] = np.dot(gradV, normal_unit)


# Log-scale transform (with small offset)
attn_weights_log = np.log10(attn_weights_np + 1e-10)

# Normalize for clarity
attn_weights_log_norm = (attn_weights_log - attn_weights_log.min()) / (attn_weights_log.max() - attn_weights_log.min())
normal_derivative_norm = (normal_derivative - normal_derivative.min()) / (normal_derivative.max() - normal_derivative.min())

# Plot clearly marking high and low points
plt.figure(figsize=(12, 5))
plt.plot(t_eval, attn_weights_log_norm, label='Log-scaled Attention Weights', color='blue')
plt.plot(t_eval, normal_derivative_norm, label='Normalized Normal Derivative of Lyapunov', color='red', linestyle='--')

# Mark the global max attention clearly
plt.scatter(t_eval[high_idx_sample], attn_weights_log_norm[high_idx_sample], color='darkred', s=100, label='Max Attention Point')

# Mark the global min attention clearly
plt.scatter(t_eval[low_idx_sample], attn_weights_log_norm[low_idx_sample], color='darkblue', s=100, label='Min Attention Point')

plt.xlabel('Time')
plt.ylabel('Normalized Value')
plt.title('Attention (Log Scale) vs. Normal Derivative of Lyapunov\n(with Max and Min Attention Points)')
plt.legend()
plt.grid(True)
plt.tight_layout()
plt.show()
\end{lstlisting}

 \bibliographystyle{plain}
 \bibliography{references}
 
\end{document}